\documentclass[11pt,twoside]{amsart} \textwidth 12 cm

\numberwithin{equation}{section}

\begin{document}

\title{On topology estimation of  submanifolds in Riemannian manifolds  by random points }
\author{R. Mirzaie}

\begin{abstract}
We show that, by sampling a sufficiently large number of random points in a neighborhood of a compact submanifold 
 $M$ of a Riemannian manifold $N$, one can recover the topology of $M$ with high confidence. 
 This holds under the assumptions on the curvatures of $M$ and $N$.

\end{abstract}
\thanks{{\scriptsize R. Mirzaie, Department of Pure Mathematics,
Faculty of Science, Imam Khomeini International University,
Qazvin, Iran \flushleft  r.mirzaei@sci.ikiu.ac.ir{\scriptsize }}}
\maketitle

\pagestyle{myheadings}

\markboth{\rightline {\scriptsize  R. Mirzaie}}
         {\leftline{Reach of submanifolds }}
\section{Introduction}

Understanding the topology of a geometric object from finite and possibly noisy data is a fundamental problem in geometry, statistics, and data analysis. In many applications ranging from manifold learning and shape recognition to the analysis of complex data sets, one assumes that the observed data points are sampled near an unknown smooth manifold $M$
 embedded in a higher-dimensional ambient space $N$. A central question is whether, and under what geometric and probabilistic conditions, one can recover the topological structure of 
$M$
 from such samples with high confidence.

A classical approach to this problem constructs a geometric complex from the data, such as the union of small geodesic balls, the Čech complex, or the Vietoris–Rips complex, and then studies its homology. The foundational results of Niyogi, Smale, and Weinberger (see [14]) established that, for compact submanifolds of Euclidean space with positive reach, the homology of such complexes coincides with that of the underlying manifold with high probability, provided the sample is sufficiently dense. Subsequent work has extended these results to various settings, including noisy samples.

In this paper, we generalize these topological inference results to the setting of compact submanifolds of Riemannian manifolds. We show that if $N$ is 
a Riemannian manifold and $M$ is a compact submanifold with positive reach in $N$, 
then, under some conditions on the curvatures of $M$ and $N$, by drawing a sufficiently large number of random points near 
$M$
 one can estimate the topology of $M$ with arbitrarily high confidence. Our analysis combines tools from Riemannian geometry and  probability
  to control the geometry of tubular neighborhoods of $M$
 and to relate random sampling properties to topological properties of $M$. 

Beyond its theoretical interest, our result provides a  toll for topological data analysis in non-Euclidean contexts, such as data lying on manifolds of non-negative curvature. It also highlights the geometric quantities, specifically the reach and curvature bounds that govern the stability of topological inference under curvature and embedding effects.

\section{Preliminaries}

In what follows, $N$ is a smooth and complete  Riemannian manifold and $M$ is a connected compact submanifold of $N$. All maps are considered to be smooth.
We will denote by $B_{\epsilon}(x)$ the open ball with the center $x$ and the radius $\epsilon>0$ in $N$. If $x \in M$, then the similar open ball
in $M$ will be denoted by $B_{\epsilon}^{M}(x)$.\\\\
The main tool in our proofs is the second variation formula which measures the change in the arch length under the small displacements.
 We refer to [15] for definitions and further results.\\\\
Let $\sigma: [0,1]
\to N$ be a geodesic segment and $\chi: [0,1] \times (-\delta,
\delta) \to N$ be a smooth map such that $\chi(t,0)=\sigma(t)$.
Then, $\chi$ is called a variation of $\sigma$. Denote by $L(s)$
the length of the curve $\sigma_{s}:[0,1] \to N$, $ t \to
\chi(t,s)$ ($\sigma_{0}=\sigma$) and let $V(t)=\frac{d}{ds}
\chi(t,s)_{_{|s=0}}=\chi_{_{s}}(t,0)$, $A(t)=\chi_{_{ss}}(t,0)$.
By the  second variation formula we have:

\[ L''(0)=\frac{1}{L(0)}\int_{0}^{1}\{<V'^{\perp},V'^{\perp}>-<R_{V\sigma'}V, \sigma'>\}dt\]\[+\frac{1}{L(0)}(<\sigma'(1), A(1)>-<\sigma'(0), A(0))>).\]
Where, $(V')^{\perp}$ is the component of $V'$ normal to
$\sigma'$.
The following lemma will play a central rule in our proofs.\\\\
{\bf Lemma 2.1.}  {\it  Let  $p, q \in N$ 
  such that $d(q,p)=\lambda>0$ and  $\sigma: [0,1] \to N$ be a  geodesic with the length $\lambda$ joining $p$ to $q$. Suppose that $\gamma: (-\delta,\delta) \to N$, $\delta>0$, is a unit speed geodesic such that
 $\gamma(0)=p$, and $d(p,q)$ is the maximum of the set $\{d(\gamma(t),q): t \in (-\delta,\delta)\}$. Then, there is a vector field $V$ along $\sigma$ such that
\[ \lambda^{2} \int_{0}^{1}\kappa(V,\sigma')(1-t)^{2}dt \geq 1.  \]  
Where $\kappa(V,\sigma')$ is the sectional curvature of $N$ along the plane generated by $V$ and $\sigma'$. }
\begin{proof}
First of all notice that since $d(p,q)$ is the maximum of the set $\{d(\gamma(t),q): t \in (-\delta,\delta)\}$, then $\gamma'(0)$ is normal to $\sigma'(0)$.  Denote by $\bar{\nabla}$ the Levi-Civita connection of $N$. 
 Let $u=\gamma'(0)$ and denote by $U(t)$  the vector field along $\sigma$ obtained by the
parallel transformation of $u$ along $\sigma$. Put $V(t)=(1-t)U(t)$ and  consider a variation $\chi(t,s)$, $(t,s) \in [0,1] \times (-\epsilon, \epsilon)$ of $\sigma$
such that $\chi_{s}(t,0)=V(t)$ and $\chi(0,s)=\gamma(s)$ ( such a variation exists for a sufficiently small $\epsilon<\delta$). Since $V(1)=0$,  then $\chi(1,s)$ is a constant point which must be $q$ and $A(1)=0$.  Also, we get from
$A(t)=\chi_{ss}(t,s)_{|s=0}$ that \[A(0)=\gamma''(0)=0 \ \ (1)\]
Denote by  $L(s)$ the length of the curve $t \to \chi(t,s)$. Thus, $L(0)$ is equal to the $\lambda$, and by the second variation formula we have
\[ L''(0)=\frac{1}{\lambda}\int_{0}^{1}\{<V'^{\perp},V'^{\perp}>-<R_{V\sigma'}V, \sigma'>\}dt.\]
Since $U$ is parallel along $\sigma$, then $\bar{\nabla}_{\sigma'}U=0$. Thus,
\[ V'(t)=\bar{\nabla}_{\sigma'}V=\bar{\nabla}_{\sigma'}(1-t)U(t)=-U(t)+(1-t)\bar{\nabla}_{\sigma'}U=-U(t) \ \ \ (2) \]
Since $U(0)(=u)$ is unit and normal to $\sigma$ at $\sigma(0)(=p)$,  $U(t)$ is unit and normal to $\sigma$ at $\sigma(t)$.
Thus,
\[ <V'^{\perp}, V'^{\perp}>=<U,U>=1 \ \ \ (3)\]
We have also
\[ <R_{V\sigma'}V, \sigma'>=\kappa(V,\sigma')<V,V><\sigma',\sigma'>=\]\[\kappa(V,\sigma')<(1-t)U,(1-t)U><\sigma',\sigma'>=\kappa(V,\sigma')(1-t)^{2} \lambda^{2} \ \ \ (4) \]
 Now, we get from (1)-(4) that
\[  L''(0)=\frac{1}{\lambda}\int_{0}^{1}(1-\kappa(V,\sigma')(1-t)^{2}\lambda^{2})dt \ \ \ (5) \]
$L(s)$ is maximum  at $s=0$, then $L''(0)\leq 0$ and we get the result from (5).

\end{proof}

{\bf Definition 2.2.}
 Let $M$ be a compact submanifold of a Riemannian manifold
$N$. The medial axes of $M$ in $N$ which we denote
 by $med(M,N)$
is the closure of the following subset of $N$.
\[\{ x \in N:   \exists p \neq p'\in  M \ such \ that \
d(x,M)=d(x,p)=d(x,p')\}\]
 The reach of $M$ in $N$ which we denote by $\tau$, is
defined as follows:
\[ \tau=inf\{ d(x,M): x \in med(M,N)\}.\]

It is easy to show that  if $reach(M,N)=\tau>0$ then  for all positive
numbers $r<\tau$, the $r$-tube around $M$, defined by $\{x \in N:
d(x,M)<r\}$ is an embedded  submanifold of $N$. If $y \in T_{r}(M)$ we will denote by $\pi(y)$ the unique nearest point of $M$ to $y$.

\section{Results}
Let $\epsilon$ be a positive number. A subset $X$ of $M$ is called $\epsilon$-dense in $M$ if for each point $y \in M$, 
$B^{M}_{\epsilon}(y,\epsilon) \cap X\neq \emptyset$. If $X \subset N$, then $X$ is called $\epsilon$-dense with respect to $M$ if
for each $y \in M$, $B_{\epsilon}(x)\neq \emptyset$\\\\
We recall that the convexity radius of $N$ is the supremum of the all numbers $r>0$ with the property that for all $x \in N$, $B_{r}(x)$
is geodesically convex.\\\\
If $X$ is a topological space and $Y \subset X$, then a countinuous map $F: X \times [0,1] \to X$ is called a contraction of $X$ to $Y$ if
\[F(x,0)=x, \ \ F(x,1) \in Y, \ \ F(y,t) \in Y, \ \ x \in X, y \in Y, t \in [0,1].   \]
Then, we say that $X$ is contractible to $Y$. It is well known that if $X$ is contractible to $Y$ then $X$ and $Y$ have the same topological properties (
have the same homology groups).\\\\
{\bf Theorem 3.1.} {\it Let  $N$ be a Riemannian manifold of convexity radius $\eta>0$, $M$ be a compact 
submanifold of $N$ such that $reach(M,N)=\tau>0$. Let $A$ be  any finite collection of points $\{x_{1},...,x_{n}\} \in M$  such that it is
$\epsilon$ dense in $M$. If $0<\epsilon <min\{\tau, \eta\}$ then $\bigcup_{i}B_{\epsilon}(x_{i})$ is contractible on $M$.}
\begin{proof}

 Since $A$ is $\epsilon$-dense in $M$, then $M \subset \bigcup_{i}B_{\epsilon}(x_{i})$. If $y \in B_{\epsilon}(x)$, $x\in A$, we get from $\epsilon < \tau$ that $\pi(y)$ is the unique nearest point of $M$ to $y$.
Let $\gamma_{y}:[0,1] \to N$ be the unit speed minimizing geodesic joining $y$ to $\pi(y)$. Since $\epsilon < \eta$, then $B_{\epsilon}(x)$ is convex. Consequently, the image of
$\gamma_{y}$ is included in $B_{\epsilon}(x)$. Put
\[ F_{x}: B_{\epsilon}(x) \times I \to B_{\epsilon}(x), \ \ F_{x}(y,t)= \gamma_{y}(t).\]
$F_{x}$ is a contraction of $B_{\epsilon}(x)$ to $B_{\epsilon}(x) \cap M$.\\
 (*) It is clear that if   $B_{\epsilon}(x_{1}) \cap B_{\epsilon}(x_{2})\neq \emptyset$, $x_{1}, x_{2} \in A$, then 
for all $y \in B_{\epsilon}(x_{1}) \cap B_{\epsilon}(x_{2})$, $F_{x_{1}}(y,t)=F_{x_{2}}(y,t)$.\\
Now, consider the function $F:\bigcup_{i}B_{\epsilon}(x_{i}) \times I \to \bigcup_{i}B_{\epsilon}(x_{i})$ defined as follows:\\
If $y \in \bigcup_{i}B_{\epsilon}(x_{i})$, consider an $x \in A$ such that $y \in B_{\epsilon}(x)$. Put $F(y,t)=F_{x}(y,t)$. By (*), $F$ is well defined and it
 is a contraction of $\bigcup_{i}B_{\epsilon}(x_{i})$ on $M$.
\end{proof}

 If $N$ is simply connected of nonpositive curvature, then $\eta=\infty$. Thus, we have the following corollary.\\
 
 {\bf Corollary 3.2.}  {\it Let $N$ be a simply connected  Riemannian manifold of nonpositive curvature and let $A$ be  any finite collection of points $\{x_{1},...,x_{n}\} \in M$  which is
$\epsilon$ dense in $M$. If $0<\epsilon <\tau$, then  $\bigcup _{x_{i} \in A}B_{\epsilon}(x_{i})$ is contractible on  $M$.} \\\\

 The above theorem provides a tool to study the topology of $M$ by using a suitable finite subset of $M$. The topology of $\bigcup _{x \in A}B_{\epsilon}(x)$ can be estimated by combinational arguments
 on the points of $A$. If  the points of $A$ are chosen  randomly, it is necessary to make sure that the union of the open balls  $B_{\epsilon}(x)$, $x \in A$, cover
 $M$.  Consequently, the estimation of the topology becomes more reliable as    the cardinality of $A$ increases. In this section we  relate the number of the points of $A$ 
 to the confidence in accurately  estimating of the topology of $M$.\\\\

 {\bf  Lemma 3.3.} {\it  Let $p \in (0,1]$ and  $\epsilon>0$, $M$ be a compact submanifold of a Riemannian manifold $N$, and $A=\{a_{1},...,a_{l}\}$ be a subset of $M$, and the scalar curvature of $M$ be bounded from the above  by $s$. Then, there is a number $\phi=\phi(p,s,\epsilon)$ such that if $l>\phi$ then with the probability at least $p$,  $A$ is $\epsilon$-dense in $M$. }
\begin{proof}

Let $dimM=m$ and consider a subset $\{y_{1},...,y_{k}\}$ of the points of $M$ such that $\bigcup_{j}B^{M}_{\frac{\epsilon}{3}}(y_{j})$ covers $M$.
$A$ is $\epsilon$-dense  if for each point $y_{j}$ there exists $a_{i} \in A$ such that $a_{i} \in B^{M}_{\frac{\epsilon}{3}}(y_{j})$.
The probability that the point $a_{i}$ belong to $B_{\frac{\epsilon}{3}}(y_{j})$ is equal to 
\[P_{i,j}=\frac {vol(B^{M}_{\frac{\epsilon}{3}}(y_{j}))}{vol(M)} \ \ \ (1)\]
Let $P'_{j}$ be the probability that there is no $a_{i}$  belonging  to $B^{M}_{\frac{\epsilon}{3}}(y_{j})$. Then,
\[ P'_{j}=(1-P_{1,j})(1-P_{2,j})...(1-P_{l,j}) \ \ \ (2)\]

Let $P'$ be the probability that there is a point $y_{j}$ such that there is no $a_{i}$ belonging to $B^{M}_{\frac{\epsilon}{3}}(y_{j})$. Then,
\[ P'\leq \sum_{j} P'_{j} \ \ \ (3)\]
Therefore, the probability $P$ that for each $j$ there exists $i$ such that $a_{i}$ belongs to $B^{M}_{\frac{\epsilon}{3}}(y_{j})$ is
\[P=1-P'\geq 1-\sum_{j} P'_{j} \ \ \ (4)\]
Now, we estimate $P'_{j}$.
We have {\bf( see []):}
\[ vol(B^{M}_{\frac{\epsilon}{3}}(y_{j}))=v_{m}(\frac{\epsilon}{3})^{2}(1-\frac{S(y_{j})}{6(m+2)})+v_{m}O(\frac{\epsilon^{4}}{3^{4}}) \ \ \ (5)\]
Where, $v_{m}$ is the volume of  standard sphere of radius one in $R^{m}$ and $S(y_{j})$ is the scalar curvature of $M$ at the point $y_{j}$. Thus,
\[ P_{i,j}\geq \frac{1}{vol(M)}v_{m}(\frac{\epsilon}{3})^{2}(1-\frac{S(y_{j})}{6(m+2)}) \ \ \ (6) \]
Since by assumption $S(y_{j})\leq s$, then
\[P_{i,j} \geq \frac{1}{vol(M)}v_{m}(\frac{\epsilon}{3})^{2}(1-\frac{s}{6(m+2)}) \ \ \ (7)\]
Then,
\[ P'_{j}\leq [1-\frac{1}{vol(M)}v_{m}(\frac{\epsilon}{3})^{2}(1-\frac{s}{6(m+2)})]^{l} \ \ \ (8)\]
Which implies that 
\[ P'\leq k[1-\frac{1}{vol(M)}v_{m}(\frac{\epsilon}{3})^{2}(1-\frac{s}{6(m+2)})]^{l} \ \ \ (9)\]
Therefore,
\[ P\geq 1-k[1-\frac{1}{vol(M)}v_{m}(\frac{\epsilon}{3})^{2}(1-\frac{s}{6(m+2)})]^{l} \ \ \ (10)\]
Let us assume that $k$ is the minimal number of balls $B_{\frac{\epsilon}{3}}(y_{j})$ covering $M$.
We find an upper bound for $k$.\\
 It is possible to choose a collection  of  balls of $M$ covering
$M$  in such a way that each point $x \in M$ belongs to at most
$m+1$ balls. Then, $\sum_{j}vol(B_{\frac{\epsilon}{3}}(y_{j})) \leq (m+1) vol(M)$.
In other way by (5) we have 
\[vol(B_{\frac{\epsilon}{3}}(y_{j}))\geq v_{m}(\frac{\epsilon}{3})^{2}(1-\frac{s}{6(m+2)}) \ \ \ (11)\]
Thus,
\[ kv_{m}(\frac{\epsilon}{3})^{2}(1-\frac{s}{6(m+2)})\leq (m+1) vol(M) \ \ \ (12)\]
Therefor,
\[ k\leq (m+1) vol(M)(v_{m}(\frac{\epsilon}{3})^{2})^{-1}(1-\frac{s}{6(m+2)})^{-1} \ \ \ (13)\]
we get from (10) and (13) that 
\[  P\geq g((m,s,l,\epsilon)\]
such that $g(m,s,l,\epsilon)$ is equal to
\[ 1-{(m+1) vol(M)(v_{m}(\frac{\epsilon}{3})^{2})^{-1}(1-\frac{s}{6(m+2)})^{-1}} \times \]\[[1-\frac{1}{vol(M)}v_{m}(\frac{\epsilon}{3})^{2}(1-\frac{s}{6(m+2)})]^{l} \]
Now, put 
\[g(m,s,l)\geq p \ \ \ (14)\] by a simple computations  we get that
(14) is equivalent to  $l\geq \phi(m,\lambda,s,\epsilon)$ for a suitable function $\phi$.

\end{proof} 
Now, putting together Theorem 3.1 and Lemma 3.3 , we get the following theorem.\\\\
{\bf Theorem 3.4.}{ \it  Let $p \in (0,1]$ and  $N$ be a Riemannian manifold of convexity radius $\eta$, $M$ be a compact 
submanifold of $N$ with the scalar curvature bounded from the above by $s$ and $reach(M,N)=\tau>0$, and let $\epsilon <min\{\tau, \eta\}$. 
Then, there is a number $\phi=\phi(p,s,\epsilon)$ such that for each finite collection of points $\{x_{1},...,x_{l}\} \in M$ , if $l>\phi$ then with the probability at leat $p$,  $\bigcup _{x \in A}B_{\epsilon}(x)$ is contractible on  $M$.}

\section{ Topology estimation from noisy data}

In Theorem 3.4, if the points of $A$ are not necessarily  on $M$ ( there is noise in selected points) it is possible again to 
estimate the topology of $M$ by $A$  under suitable conditions on the selected points. The following theorem 
describes a situation  where we can obtain topological properties of $M$ from a noisy subset $A$ by high confidence. We suppose that the selected 
points are spread around $M$ (i.e, they belong to an $r$-tube around $M$ with suitably selected positive number $r$). As the previous theorems, the reach of $M$ in $N$
plays an important rule.\\\\
We will use the following lemma.\\

 {\bf Lemma 4.1.} {\it  Let $p \in (0,1]$ and $\epsilon>0$, and $M$ be a compact submanifold of a Riemannian manifold $N$, and  $A=\{a_{1},...,a_{l}\}$ be a subset of $N$ included in $r$-tube around $M$, $r>0$. Suppose that the scalar curvature of $N$ is bounded from the above  by $s$. Then, there is a number $\phi=\phi(r,p,s,\epsilon)$ such that if $l>\phi$ then with the probability at least $p$,  $A$ is $\epsilon$-dense with respect to $M$. }
\begin{proof}
The proof is similar to the proof of Lemma 3.3, with a little changes as follows:\\
 1. We replace  $B^{M}_{\epsilon}(y_{i})$ by $B_{\epsilon}(y_{i})$ ( open ball in $N$). \\
 2. $vol(M)$ must be replaced by $vol(T_{r}(M))$, where  $T_{r}(M)$ is the $r$-tube around $M$.\\
 3. $S(y_{j})$ is the scalar curvature of $N$.\\
 4. $m$ the dimension of $M$  must be replaced by $n$ the dimension of $N$.\\
 5. $v_{m}$ must be replaced by $v_{n}$ ( the volume of standard sphere of radius one in $R^{n}$).\\

\end{proof}

{\bf Lemma 4.2.} {\it Let $M$ be a compact submanifold of a Riemannian manifold $N$ such that $\tau=reach(M,N)>0$. Let $0<\epsilon <\frac{\tau}{2}$ and 
$C=\{x_{1},...,x_{l}\}$ be a subset of $N$ such that:\\
(a) $C \subset T_{\frac{\tau}{2}}(M)$.\\
(b) $C$ be $\frac{\epsilon}{2}$-dense with respect to $M$ (i.e, for each $y \in M$, $B_{\frac{\epsilon}{2}}(y) \cap C\neq \emptyset)$. \\
(c) The sectional curvature of $N$ be upper bounded by $  \frac{1}{25 \tau^{2}}$.\\
Then, $\bigcup_{i}B_{\epsilon}(x_{i})$ is contractible on $M$.} \\

\begin{proof}
Let $W=\bigcup_{i}B_{\epsilon}(x_{i})$. We show that $W$ is contractible on $M$.
Since $C$ is $\frac{\epsilon}{2}$-dense with respect to $M$, then $M \subset W$. 
Since $\epsilon <\frac{\tau}{2}$, by (a) and  definition of the reach, for each point $z \in W$  there is a unique point $\hat{z} \in M$ such that
$d(z,M)=d(z,\hat{z})$. let $\gamma=\gamma_{z}:[-1,1] \to N$, be the unique  geodesic in $N$ joining $z$ to $\hat{z}$.
We show that for all $t \in [-1,1]$, $\gamma_{z}(t)$ belongs to $W$.\\
If there is a point $z'=\gamma(s)$, $s \in [-1,1]$, such that $z' \notin W$ then we ge a contradiction as follows:\\

Let $\Gamma=\{\gamma(s): s \in [-1,1], \gamma(s) \notin W\}$.  
Choose  a $z'=\gamma(s_{0}) \in \Gamma$ with the property  that 
\[ d(z',\overline{W})=max\{ d(\gamma(t),\overline{W}), \ \ t \in [-1,1]\} \ \  (1)\]
$\overline{W}$ is the closure of $W$ in $N$ which  is compact and connected. By the compactness of $\overline{W}$ and $[-1,1]$, such a $z'$ exists which without lose 
of the generality we assume that $z'=\gamma(0)$. 
Consider a point $x' \in \overline{W}$ such that 
\[ d(z',x')=d(z', \overline{W})=inf\{d(z',x): x \in \overline{W }\} \ \ (2)\]
 It is easy to show that by (1), there is small number $0<\delta<1$, such that $d(z',x')$ is maximum of $d(\gamma(t),x')$ for $t \in[-\delta,\delta]$. Let $\sigma:[0,1] \to N$
 be the geodesic with the length $\lambda=d(x',z')$ joining $x'$ to $z' $.
By Lemma 2.1, we have:
\[ \lambda^{2} \int_{0}^{1}\kappa(V,\sigma')(1-t)^{2}dt \geq 1 \ \ (3)  \]
Since $C$ is $\frac{\epsilon}{2}$-dense with respect to $M$, there is a point $x'' \in C \subset W$ such that 
\[ d(\hat{z},x'')< \frac{\epsilon}{2} \ \ (4)\]
By (2) we have 
\[ d(z',x')\leq d(z',x'')\leq d(z',\hat{z})+d(\hat{z},x'')  \ \ (5) \]
By definition of the reach, (a) and the assumption $\epsilon <\frac{\tau}{2}$, we get that $d(z,\hat{z})<\tau$. Thus
\[ d(z',\hat{z})=d(\gamma(0),\hat{z})\leq d(\gamma(-1),\hat{z})=d(z,\hat{z})<\tau \ \ (6) \]
(4),(5), (6) yield to 
\[ \lambda=d(z',x') \leq \tau+\frac{\epsilon}{2}  \]
Since by assumption $\epsilon < \frac{\tau}{2}$, then
\[\lambda < \frac{5}{4} \tau \ \ (7)\]

We get from (3) and (7) that
\[\int_{0}^{1}\kappa(V,\sigma')(1-t^{2})dt \geq \frac{16}{25 \tau^{2}} \ \ (8)\]
In other way, by the assumption (c) of the theorem we have:
\[\int_{0}^{1}\kappa(V,\sigma')(1-t)^{2}dt< \frac{1}{25\tau^{2}}\int_{0}^{1}(1-t)^{2}dt= \frac{1}{25\tau^{2}} \times\frac{2}{3} \ \ (9)\]
By (8) and (9) we have  a contradiction.\\
At the end, the following map gives a contraction of $W$ on $M$.
\[ F: W \times [-1,1] \to W, \ \ F(z,t)=\gamma_{z}(t).\]

\end{proof}
 
Now, we get the following theorem from Lemma 4.1 and Lemma 4.2.\\\\
{ \bf Theorem 4.3.} {\it  Let $p \in (0,1]$ and $M$ be a compact submanifold of a Riemannian manifold $N$ such that $\tau=reach(M,N)>0$, and  let  the sectional curvature of $N$ be upper bounded by $  \frac{1}{25 \tau^{2}}$, and let  $\epsilon <\frac{\tau}{2}$. Then,
there is a number $\phi=\phi(p,\eta,\epsilon)$ such that for each collection 
$C=\{x_{1},...,x_{l}\}$ of the points in $N$ selected from $T_{\frac{\tau}{2}}(M)$, if  $l>\phi$
then with the probability at least $p$, $\bigcup_{i}B_{\epsilon}(x_{i})$ is contractible to  $M$.  }
\begin{proof}
If we add the condition that the   scalar curvature $S$ of $N$ be upper bounded by $s$, then  by Lemma 4.1 and Lemma 4.2, there is a number
$\phi=\phi(p,s,\epsilon)$ with the desired property. But, we have $S\leq n(n-1)K_{max}$, where $K_{max}$ is the maximum of the sectional curvatures on $M$.
Thus, the upper bound $\frac{1}{25 \tau^{2}}$ for sectional curvatures imposes an upper bound for the scalar curvature.  Then, $\phi$ can be written as a function
of $(p, \eta,\epsilon)$ (i.e, $\phi=\phi(p,\eta,\epsilon)$).
\end{proof}
{\bf  Remak 4.4.} By Theorems 3.4 and 4.3, the topological properties of a submanifold 
$M$ of $N$ can be recovered with arbitrarily high confidence from a finite point cloud
$\{x_{1},...,x_{l}\}$ consisting either of points sampled on $M$ or points sampled near $M$. To analyze the topology of $\bigcup_{i}B_{\epsilon}(x_{i})$
 in practice, one constructs a filtration of simplicial complexes whose vertex set is the sample points 
$\{x_{1},...,x_{l}\}$. Examples include the Čech complex, Vietoris–Rips complex, or witness complexes.
 Theorems 3.4 and 4.3 guarantee that, under sampling and geometric conditions controlled by the reach $\tau$, the topology of these complexes recovers the topology of $M$.\\
{\bf Role of the Ambient Curvature}.\\
{\bf Euclidean ambient space $R^{n}$}. The ambient space has zero sectional and scalar curvature. In this case, the results reduce to classical theorems such as those of Niyogi–Smale–Weinberger (2008) and related work in topological inference.\\
{\bf Ambient manifold $N$ with non-positive curvature}.
In Theorem 4.3, the assumption
 on the sectional curvature of $N$ is automatically satisfied. Therefore, the theorem holds for every compact submanifold $M \subset N$ with positive reach.\\
{\bf Ambient manifold $N$ with positive sectional curvature}.\\
 In this case, the upper bound $\kappa_{N}\leq \frac{1}{25 \tau^{2}}$ 
 becomes essential. This ensures that the ambient curvature does not distort geodesics of $M$ too much relative to the reach of $M$.

\end{document}